\documentclass[12pt,a4paper]{scrartcl} 
\usepackage[utf8]{inputenc}
\usepackage[english,russian]{babel}
\usepackage{indentfirst}
\usepackage{misccorr}
\usepackage{graphicx}
\usepackage{amsmath}
\begin{document}
	\begin{center}
		\large
	    \textbf{PRIME COUNTING FUNCTION}
	    \\
	    IGOR TURKANOV   
	\end{center}
\medskip
\begin{quote}
	\large
	ABSTRACT
\\	
The theorem below gives another way of computing the distribution prime counting function without using recursion and the values of Prime numbers
\\	
\end{quote}
\large
\chapter{THEOREM}
\\
\\
The prime counting function is the function  $\pi(n)$, giving the number of primes less than or equal to a given number $n$, in explicit form is expressed by the formula:

\begin{displaymath}
 \pi(n)=n-1-\sum\limits_{i=2}^{\left[\sqrt{n}\right]}\left(\left[\frac{n}{i}\right]-i+1\right)+
\end{displaymath}
\begin{displaymath}\label{eq:solv}
 \sum\limits_{s=2}^{\left[\sqrt{n}\right]}(-1)^{s}
 \sum\limits_{1<i_1<i_2<\dotsc<i_s\leqslant\left[\sqrt{n}\right]}
 \left(\left[\frac{n}{LCM(i_1,i_2,\dotsc,i_s)}\right]-\left[\frac{i_s^2-1}
 {LCM(i_1,i_2,\dotsc,i_s)}\right]\right)  
\end{displaymath}

Where the $[x]$ is the floor function of $x$,
$LCM(i_1,i_2,\dotsc,i_s)$ the least common multiple of positive integers $i_1,i_2,\dotsc,i_s$.
\\
\\
\chapter{PROOF:}	
\\
\\
Function $\pi(n)$ is equal to a difference  $n-1$   ($1$ - by definition not prime) and numbers of compound numbers, each of which is $i\cdot j$ with conditions: $i\cdot j \leqslant n,\;\; i,j \in \rm N_1$, where through $\rm N_1$ we shall designate set of natural numbers, greater  $1$.

Let's add a condition $j \geqslant i $ which does not limit a generality.

Thus, all natural compound numbers, smaller or equal $n $, form set:
\begin{center}
$X = \bigl\{x\;\lvert\; x\leqslant n,\; x=i\cdot j,\; j\geqslant i,\;\; x,i,j \in \rm N_1\bigr\}$
.\end{center}

By definition $\pi(n)=n-1-\lvert X \rvert$.

Let's designate through $X_i $ set of natural compound numbers of a type $i \cdot j $, not surpassing $n $, with fixed $i \in N_1 $:
\begin{center}
$X_i = \bigl\{x\;\lvert\; x\leqslant n,\; x=k\cdot j,\;k=i,\; j\geqslant k,\;\; x,j,k \in \rm N_1\bigr\}$.
\end{center}

Note that $X_i=\varnothing, \; \; \forall i > \left[\sqrt{n} \right]$.
\\
\\
For clarity, draw a table, for example, for $n=11$:
\begin{center}
$\begin{smallmatrix}j\\11 && 11\\10 && 10\\9 && 9\\8 && 8\\7 && 7\\6 && 6\\5 && 5 && 10\\4 && 4 && 8\\3 && 3 && 6 && \mathbf{9}\\
2 && 2 && \mathbf{4} && 6 && 8 && 10\\
1 && \mathbf{1} && 2 && 3 && 4 && 5 && 6 && 7 && 8 && 9 && 10 && 11\\ \; && 1 && 2 && 3 && 4 && 5 && 6 && 7 && 8 && 9 && 10 && 11 && i \end{smallmatrix}$
\end{center}

In an example $X_2 = \bigl\{4,6,8,10 \bigr\},\;\;X_3 = \bigl\{9 \bigr\}$.
\\
\\
In the table $X_k $ is all numbers from $k $-th column, it is more or equal $k^2 $.
\\
It's clear that:

\begin{center}
	$X=X_2\cup X_3\cup X_4\cup \dotsb \cup X_{\left[\sqrt{n}\right]}$.
\end{center}

In columns numbers $i_1 $ and $i_2 $ there can be identical numbers, i.e. $X_{i_1}\cap X_{i_2}$ is not always empty.
\\
In accordance with the Inclusion-Exclusion Principle~\cite{Andrews-1971}:
\begin{center}
\[\lvert X\rvert=\sum\limits_{i=2}^{\left[\sqrt{n}\right]} \lvert  X_i\rvert-\sum\limits_{s=2}^{\left[\sqrt{n}\right]}(-1)^{s}\sum\limits_{1<i_1<i_2<\dotsc<i_s\leqslant\left[\sqrt{n}\right]}\lvert X_{i_1} \cap X_{i_2} \cap \dotsb \cap X_{i_s}\rvert\]
\end{center}

Finding power sets, the right-hand side of this equation,
we will find $\pi(n)$.
\bigskip

Let us prove the following statement.
\pagebreak
\\
\\
\chapter{STATEMENT}
\\
\\
For $\forall i,s,i_1,i_2,\dotsc,i_s \in N_1,\;1<i_1<i_2<\dotsc<i_s\leqslant\left[\sqrt{n}\right]$,\\
$X_i = \bigl\{x\;\lvert\; x\leqslant n,\; x=k\cdot j,\;k=i,\; j\geqslant k,\;\; x,j,k \in \rm N_1\bigr\}$  
right:
\begin{center}
$\lvert X_{i_1} \cap X_{i_2} \cap \dotsb \cap X_{i_s}\rvert=
\left[\frac{n}{LCM(i_1,i_2,\dotsc,i_s)}\right]-\left[\frac{i_s^2-1}{LCM(i_1,i_2,\dotsc,i_s)}\right],
$
\end{center}
Where the $[x]$ is the floor function of $x$,
$LCM(i_1,i_2,\dotsc,i_s)$ the least common multiple of positive integers $i_1,i_2,\dotsc,i_s$.
\\
\\
\chapter{PROOF:}
\\
\\
Let $Y_i(n)$ the set of natural composite numbers of the form $i\cdot j$, not exceeding $n$, with a fixed $i \in N_1$, without the condition $j\geqslant i$, $j \in N$:
\begin{center}
	$Y_i(n) = \bigl\{y\;\lvert\; Y\leqslant n,\; y=k\cdot j,\;k=i,\;\; y,k \in N_1, \; j \in N\bigr\}$.
\end{center}
$Y_k(n)$ - set of numbers, standing in the $k$-th column of the table in question. 
In an example: 
$Y_2(11) = \bigl\{2,4,6,8,10 \bigr\},\;\;Y_3(11) = \bigl\{3,6,9 \bigr\},
\;\;Y_2(3) = \bigl\{2\bigr\},\;\;Y_3(8) = \bigl\{3,6\bigr\}$.
\\
From the above definition that sets $\forall s,n,i_1,i_2,\dotsc,i_s \in N_1$, $1<i_1<i_2<\dotsc<i_s\leqslant\left[\sqrt{n}\right]$ true equality:
\begin{displaymath}    
    \lvert X_{i_1} \cap X_{i_2} \cap \dotsb \cap X_{i_s}\rvert=
\end{displaymath}
\begin{displaymath}   
    \lvert Y_{i_1}(n) \cap Y_{i_2}(n) \cap \dotsb \cap Y_{i_s}(n)\rvert -
    \lvert Y_{i_1}(i_s^2-1) \cap Y_{i_2}(i_s^2-1) \cap \dotsb \cap Y_{i_s}(i_s^2-1)\rvert   
\end{displaymath}
In the left part of equality there is a quantity of numbers in a column with the maximal index $i_s $, conterminous with numbers from columns with numbers $i_1, i_2, \dotsc, i _ {s-1} $, a type $i_s \cdot j $, not surpassing $n $, with a condition $j \geqslant i $.
\\
In the right part – a difference of quantity of the same numbers in the same $i_s $-th column without a condition $j \geqslant i_s $ and quantities of the same numbers which size less or is equal $i_s^2-1 $, that from a way of construction of the table to equivalently condition $j <i_s $.

We prove by induction on the index of $s$, that $\forall s, n, i_1, i_2, \dotsc, i_s \in N_1 $:
\begin{displaymath}
Y_{i_1}(n) \cap Y_{i_2}(n) \cap \dotsc \cap Y_{i_s}(n) =\\
\end{displaymath}
\begin{displaymath}
\bigl\{y\;\lvert\; y\leqslant n,\; 
\exists\; m \in N:\; y=LCM(i_1,i_2, \dotsc ,i_s) \cdot m,\;\; y \in N_1\bigr\}
\end{displaymath}

and then:
\begin{displaymath}
\lvert Y_{i_1}(n) \cap Y_{i_2}(n) \cap \dotsb \cap Y_{i_s}(n)\rvert =
\left[\frac{n}{LCM(i_1,i_2,\dotsc,i_s)}\right].
\end{displaymath}
\\
For $s=1$:
\\
For $\forall i \in N_1$ 
\begin{displaymath}
Y_i(n) = \bigl\{y\;\lvert\; y\leqslant n,\; y=k\cdot j,\;k=i,\;\; y,k \in N_1, \;j \in N\bigr\},
\end{displaymath} 
which is equivalent to:
\begin{displaymath}
Y_i(n) =
\bigl\{y\;\lvert\; y\leqslant n,\; 
\exists\; j \in N:\; y=LCM(i) \cdot j,\;\; y \in N_1\bigr\}
\end{displaymath}
and 
\begin{displaymath}
\lvert Y_i(n)\rvert = \left[\frac{n}{i}\right] = \left[\frac{n}{LCM(i)}\right].
\end{displaymath}
The first step of the induction is confirmed by these equations.
\\
\\
Suppose that for $s=k-1$ is true equality, we show that for $s=k$ is also true.
\begin{displaymath}
 Y_{i_1}(n) \cap Y_{i_2}(n) \cap \dotsc \cap Y_{i_{k-1}}(n) =\\
\end{displaymath}
\begin{displaymath}
\bigl\{y\;\lvert\; y\leqslant n,\; 
\exists\; m \in N:\; y=LCM(i_1,i_2, \dotsc ,i_{k-1}) \cdot m,\;\; y \in N_1\bigr\}, \end{displaymath}
\\
\begin{displaymath}
Y_{i_1}(n) \cap Y_{i_2}(n) \cap \dotsc \cap Y_{i_k}(n) =\\
\end{displaymath}
\begin{displaymath}
\bigl\{y\;\lvert\; y\leqslant n,\; 
y=i_1\cdot j_1=i_2\cdot j_2=\dotsc =i_k\cdot j_k,\;\; y \in N_1,\;j_1,j_2,\dotsc, j_k \in N_1\bigr\}=
\end{displaymath}
\begin{displaymath}
\bigl\{y\;\lvert\; y\leqslant n,\; \exists m,r \in N :\; y=LCM(i_1,i_2, \dotsc ,i_{k-1}) \cdot m = 
i_k \cdot r,\;\; y\in N_1\bigr\}=
\\
\end{displaymath}
\begin{displaymath}
\bigl\{y\;\lvert\; y\leqslant n,\; 
\exists\; m \in N:\; y=LCM(LCM(i_1,i_2, \dotsc ,i_{k-1}),i_{k}) \cdot m,\;\; y \in N_1\bigr\},
\end{displaymath}
i.e.
\begin{displaymath}
\lvert Y_{i_1}(n) \cap Y_{i_2}(n) \cap \dotsc \cap Y_{i_{k}}(n)\rvert=\\
\end{displaymath}
\begin{displaymath}
\left[\frac{n}{LCM(LCM(i_1,i_2,\dotsc,i_{k-1}),i_{k})}\right] =
\left[\frac{n}{LCM(i_1,i_2,\dotsc,i_{k})}\right].
\end{displaymath}
\\
Therefore, by induction to $\forall s \in N_1$:
\begin{displaymath}
\lvert Y_{i_1}(n) \cap Y_{i_2}(n) \cap \dotsc \cap Y_{i_{s}}(n)\rvert=
\left[\frac{n}{LCM(i_1,i_2,\dotsc,i_{s})}\right].
\end{displaymath}

As a result, 
$\lvert X_{i_1} \cap X_{i_2} \cap \dotsb \cap X_{i_s}\rvert=
\left[\frac{n}{LCM(i_1,i_2,\dotsc,i_s)}\right]-\left[\frac{i_s^2-1}{LCM(i_1,i_2,\dotsc,i_s)}\right]$, 
as proves the Statement.

Having substituted expression for capacity of intersection of sets $X_{i_1}, X_{i_2}, \dotsc, X_{i_s} $ in the formula of inclusion-exception, we shall receive the proof of the basic theorem, in view of that in the first sum for presentation $\lvert X_i\rvert =  \left[\frac{n}{LCM(i)}\right]-\left[\frac{i^2-1}{LCM(i)}\right]=
\left[\frac{n}{i}\right]-\left[\frac{i^2-1}{i}\right]=
\left[\frac{n}{i}\right]-\left[i-\frac{1}{i}\right]=
\left[\frac{n}{i}\right]-i+1$. 
\\

\renewcommand{\refname}{REFERENCES:}

\end{document}